# An analytical closed-form solution for free vibration of stepped circular/annular Mindlin functionally graded plate

M. Derakhshani [a*], Sh. Hosseini-Hashemi [b], M. Fadaee [c]

[a] *M.Sc. Student of Mechanical Engineering, School of Mechanical Engineering, Iran University of Science and Technology, Tehran, Iran.*

[b] *Professor of Mechanical Engineering, School of Mechanical Engineering, Iran University of Science and Technology, Tehran, Iran.*

[c] *Assistant Professor of Mechanical Engineering, School of Mechanical Engineering, Qom University of Technology, Qom, Iran.*

*\* m.derakhshany87@gmail.com*

## Abstract

An exact solution based on a unique procedure is presented for free vibration of stepped circular and annular functionally graded (FG) plates via first-order shear deformation plate theory of Mindlin. A power-law distribution of the volume fraction of the components is considered for the Young's Modulus and Poisson's ratio of the studied FG plate. Free vibration of the plate is solved by introducing some new potential functions and the use of separation of variables method. Finally, several comparisons of the developed model were presented with the FEA analysis, to demonstrate the accuracy of the proposed exact procedure. The effect of the geometrical parameters such as step thickness ratios and step locations on the natural frequencies of FG plates is also investigated.

**Keywords**: Free vibration; Stepped circular plate; Functionally graded material; Mindlin theory.

## 1. Introduction

Numerous studies have been performed on free vibration of circular and annular plates with variable thickness. A summary of works on free vibration of thin circular and annular plates with variable thickness was done by Leissa, most of which were developed based on the classical plate theory (CPT). One of the problems associated with using this theory is that it is not able to model the expected behavior for moderately thick plates. To resolve such problem, some models have been developed by researchers based on the theory of moderately thick and thick plates including the first shear deformation theory (FSDT), the third order shear deformation theory (TSDT) and the 3D elasticity, in which the effects of transverse shear deformation and rotary inertia are considered [1-4].

In recent years, functionally graded materials (FGM) have gained so many attentions as special composites with material properties that vary continuously through their thickness. This continuous change in material properties results in eliminating discontinuities of stresses, high resistance to temperature gradients, reduction in residual and thermal stresses, high wear resistance, and an increase in strength to weight ratio. Circular and annular FG plates have been reviewed by many researchers due to the multiple applications they can be used for. Efraim and Eisenberger, Tajeddini, Hosseini-Hashemi, and Gupta are the most famous researchers who have studied the vibration of FG circular/annular plates [5-8]. It has been found that there are very limited studies performed on the circular/annular FG plates with variables thickness and to the best of the authors' knowledge, there is no study for exact closed-form solutions of vibration analysis of stepped thickness circular/annular Mindlin FG plate.

In this study, an exact solution is presented for free vibration of stepped thickness moderately thick circular and annular FG plates. By use of some new auxiliary and potential functions along with a unique theoretical approach [7], the equations of motion are exactly solved. The accuracy of the current analytical approach is validated by comparing the obtained results of the proposed model with finite element analysis (FEA). Furthermore, the effect of the various parameters such as step thickness ratios, and step locations on the natural frequencies of FG plates is evaluated. More details about the current study could be found in [9].

## 2. Mathematical formulation

Consider an annular functionally graded plate of radius $r_n$ and thickness $h_n$ consists of $n$ steps in which the radius and thickness of $i^{\text{th}}$ annular segment is $r_i$ and $h_i$ respectively, as shown in **Figure 1**. The plate geometry and dimensions are defined in an orthogonal cylindrical coordinate system $(r,\theta,z)$. In this study, the properties of the plate are assumed to vary through the plate thickness with a power-law distribution of the volume fractions of two materials. The top surface of the first segment ($z = h_1/2$), which assumed here the thickest part of the plate, is metal-rich while the bottom surface of the same segment ($z = -h_1/2$) is ceramic-rich. By considering this assumption, Young's modulus and mass density are vary through the plate thickness by the following relations

$$E(z) = (E_m - E_c)V_f(z) + E_c \qquad \rho(z) = (\rho_m - \rho_c)V_f(z) + \rho_c \qquad (1a\text{-}b)$$

where

$$V_f(z) = \left(\frac{z}{h_1} + \frac{1}{2}\right)^g \qquad (2)$$

in which the subscripts $m$ and $c$ represent the metallic and ceramic constituents, respectively. $V_f$ shows the volume fraction and $g$ is the power-law index which takes only non-negative values. Poisson's ratio $v$ is taken as 0.3 throughout the analysis.

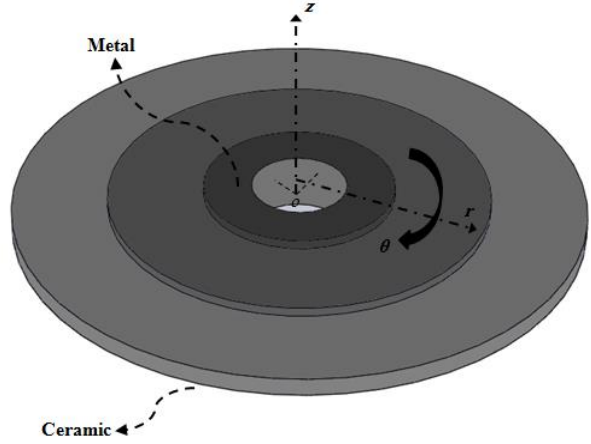

**Figure 1**. Geometry of a stepped annular FG plate

A stepped annular Mindlin plate with $n$ steps can be divided to $n$ annular plates. An annular Mindlin plate of radius $r_i$ and thickness $h_i$, in which $i$ refers to the $i^{th}$ segment of the stepped annular FG plate as shown in **Figure 1**, is considered in this model. According to the FSDT, the displacement field is used for the $i^{th}$ segment of the stepped annular Mindlin FG plate as follows

$$u^i(r,\theta,z,t) = u_0^i(r,\theta,t) + z\psi_r^i(r,\theta,t) \qquad (3a)$$

$$v^i(r,\theta,z,t) = v_0^i(r,\theta,t) + z\psi_\theta^i(r,\theta,t) \qquad (3b)$$

$$w^i(r,\theta,z,t) = w_0^i(r,\theta,t) = w^i(r,\theta,t) \qquad (3c)$$

where $u^i$, $v^i$ and $w^i$ denote the displacements in $r$, $\theta$ and $z$ directions, respectively. $u_0^i$ and $v_0^i$ denote the in-plane displacements of mid-plane in radial and circumferential directions. Moreover, $\psi_r^i$ and $\psi_\theta^i$ show the slope rotations in $r$-$z$ and $\theta$-$z$ planes at $z=0$ for $i^{th}$ segment of the plate, respectively.

The exact free vibration of circular/annular Mindlin FG plate has been studied by the first author [7], recently. In this study, the same analytical approach is used to solve the free vibration of circular/annular moderately thick FG plates. For generality and convenience in deriving mathematical formulations, the following non-dimensional terms are defined:

$$R = \frac{r}{r_n}, \quad Z_i = \frac{z}{h_i}, \quad \delta_i = \frac{h_i}{r_n}, \quad \tau_i = \frac{h_i}{h_n} \tag{4a-d}$$

For harmonic motion, the displacement fields in dimensionless forms are taken as

$$\bar{u}^i(R,\theta,Z_i) = \frac{u^i(r,\theta,z,t)}{h_i} e^{-j\omega t}, \quad \bar{v}^i(R,\theta,Z_i) = \frac{v^i(r,\theta,z,t)}{h_i} e^{-j\omega t} \tag{5a, b}$$

$$\bar{w}^i(R,\theta) = \frac{w^i(r,\theta,t)}{r_n} e^{-j\omega t}, \quad \bar{\psi}_r^i = \psi_r^i e^{-j\omega t}, \quad \bar{\psi}_\theta^i = \psi_\theta^i e^{-j\omega t} \tag{5c-e}$$

where

$$\bar{u}^i(R,\theta,Z_i) = \bar{u}_0^i + Z\bar{\psi}_r^i, \quad \bar{v}^i(R,\theta,Z_i) = \bar{v}_0^i + Z\bar{\psi}_\theta^i, \quad \bar{w}^i(R,\theta) = \bar{w}^i \tag{6a-c}$$

By substituting the stress-strain relations in polar coordinate [4] into stress resultants, which are obtained by the first author in [7] based on the FSDT, stress resultants in dimensionless forms for $i^{th}$ segment of the stepped circular FG plate are obtained as

$$\bar{N}_k^i = \frac{N_k^i}{E_c h_i} e^{-j\omega t}, \quad \bar{M}_k^i = \frac{M_k^i}{E_c h_i^2} e^{-j\omega t}, \quad \bar{Q}_k^i = \frac{Q_k^i}{E_c h_i} e^{-j\omega t} \quad k = r,\theta,r\theta \tag{7a-c}$$

in which

$$(N_r^i, N_\theta^i, N_{r\theta}^i) = \int_{-h_i/2}^{h_i/2} (\sigma_{rr}^i, \sigma_{\theta\theta}^i, \sigma_{r\theta}^i) dz \tag{8a}$$

$$(M_r^i, M_\theta^i, M_{r\theta}^i) = \int_{-h_i/2}^{h_i/2} (\sigma_{rr}^i, \sigma_{\theta\theta}^i, \sigma_{r\theta}^i) z \, dz \tag{8b}$$

$$(Q_r^i, Q_\theta^i) = \int_{-h_i/2}^{h_i/2} (\sigma_{rz}^i, \sigma_{\theta z}^i) dz \tag{8c}$$

By substituting Eqs. (7a-c) into equations of motion obtained in [7], the final dimensionless form of equations of motion are obtained as follows

$$\bar{K}_1^i \left[ \left( \frac{\partial^2 \bar{u}_0^i}{\partial R^2} + \frac{\partial \bar{u}_0^i}{R \partial R} - \frac{\bar{u}_0^i}{R^2} + \frac{\partial^2 \bar{v}_0^i}{R \partial R \partial \theta} - \frac{\partial \bar{v}_0^i}{R^2 \partial \theta} \right) + \frac{1-\nu}{2} \left( \frac{\partial^2 \bar{u}_0^i}{R^2 \partial \theta^2} - \frac{\partial \bar{v}_0^i}{R^2 \partial \theta} - \frac{\partial \bar{v}_0^i}{R \partial R \partial \theta} \right) \right]$$

$$+ \bar{K}_2^i \left[ \left( \frac{\partial^2 \bar{\psi}_r^i}{\partial R^2} + \frac{\partial \bar{\psi}_r^i}{R \partial R} - \frac{\bar{\psi}_r^i}{R^2} + \frac{\partial^2 \bar{\psi}_\theta^i}{R \partial R \partial \theta} - \frac{\partial \bar{\psi}_\theta^i}{R^2 \partial \theta} \right) + \frac{1-\nu}{2} \left( \frac{\partial^2 \bar{\psi}_r^i}{R^2 \partial \theta^2} - \frac{\partial \bar{\psi}_\theta^i}{R^2 \partial \theta} - \frac{\partial^2 \bar{\psi}_\theta^i}{R^2 \partial R \partial \theta} \right) \right] \tag{9a}$$

$$= -S_1^i \lambda_i^2 (\bar{I}_1^i \bar{u}_0^i + \bar{I}_2^i \bar{\psi}_r^i)$$

$$\bar{K}_1^i \left[ \frac{\partial \bar{u}_0^i}{R^2 \partial \theta} + \frac{\partial^2 \bar{u}_0^i}{R \partial R \partial \theta} + \frac{\partial^2 \bar{v}_0^i}{R^2 \partial \theta^2} + \frac{1-v}{2} \left( -\frac{\partial^2 \bar{u}_0^i}{R \partial R \partial \theta} + \frac{\partial \bar{u}_0^i}{R^2 \partial \theta} - \frac{\bar{v}_0^i}{R^2} + \frac{\partial \bar{v}_0^i}{R \partial R} + \frac{\partial^2 \bar{v}_0^i}{\partial R^2} \right) \right]$$

$$+\bar{K}_2^i \left[ \frac{\partial \bar{\psi}_r^i}{R^2 \partial \theta} + \frac{\partial^2 \bar{\psi}_r^i}{R \partial R \partial \theta} + \frac{\partial^2 \bar{\psi}_\theta^i}{R^2 \partial \theta^2} + \frac{1-v}{2} \left( -\frac{\partial^2 \bar{\psi}_r^i}{R \partial R \partial \theta} + \frac{\partial \bar{\psi}_r^i}{R^2 \partial \theta} - \frac{\bar{\psi}_\theta^i}{R^2} + \frac{\partial \bar{\psi}_\theta^i}{R \partial R} + \frac{\partial^2 \bar{\psi}_\theta^i}{\partial R^2} \right) \right]$$

$$= -S_1^i \lambda_i^2 (\bar{I}_1^i \bar{v}_0^i + \bar{I}_2^i \bar{\psi}_\theta^i)$$

(9b)

$$\bar{K}_2^i \left[ \left( \frac{\partial^2 \bar{u}_0^i}{\partial R^2} + \frac{\partial \bar{u}_0^i}{R \partial R} - \frac{\bar{u}_0^i}{R^2} + \frac{\partial^2 \bar{v}_0^i}{R \partial R \partial \theta} - \frac{\partial \bar{v}_0^i}{R^2 \partial \theta} \right) + \frac{1-v}{2} \left( \frac{\partial^2 \bar{u}_0^i}{R^2 \partial \theta^2} - \frac{\partial \bar{v}_0^i}{R^2 \partial \theta} - \frac{\partial \bar{v}_0^i}{R \partial R \partial \theta} \right) \right]$$

$$+\bar{K}_3^i \left[ \left( \frac{\partial^2 \bar{\psi}_r^i}{\partial R^2} + \frac{\partial \bar{\psi}_r^i}{R \partial R} - \frac{\bar{\psi}_r^i}{R^2} + \frac{\partial^2 \bar{\psi}_\theta^i}{R \partial R \partial \theta} - \frac{\partial \bar{\psi}_\theta^i}{R^2 \partial \theta} \right) + \frac{1-v}{2} \left( \frac{\partial^2 \bar{\psi}_r^i}{R^2 \partial \theta^2} - \frac{\partial \bar{\psi}_\theta^i}{R^2 \partial \theta} - \frac{\partial^2 \bar{\psi}_\theta^i}{R^2 \partial R \partial \theta} \right) \right]$$

$$-S_2^i (\bar{\psi}_r^i + \frac{\partial \bar{w}^i}{\partial R}) = -S_1^i \lambda_i^2 (\bar{I}_2^i \bar{u}_0^i + \bar{I}_3^i \bar{\psi}_r^i)$$

(9c)

$$\bar{K}_2^i \left[ \frac{\partial \bar{u}_0^i}{R^2 \partial \theta} + \frac{\partial^2 \bar{u}_0^i}{R \partial R \partial \theta} + \frac{\partial^2 \bar{v}_0^i}{R^2 \partial \theta^2} + \frac{1-v}{2} \left( -\frac{\partial^2 \bar{u}_0^i}{R \partial R \partial \theta} + \frac{\partial \bar{u}_0^i}{R^2 \partial \theta} - \frac{\bar{v}_0^i}{R^2} + \frac{\partial \bar{v}_0^i}{R \partial R} + \frac{\partial^2 \bar{v}_0^i}{\partial R^2} \right) \right]$$

$$+\bar{K}_3^i \left[ \frac{\partial \bar{\psi}_r^i}{R^2 \partial \theta} + \frac{\partial^2 \bar{\psi}_r^i}{R \partial R \partial \theta} + \frac{\partial^2 \bar{\psi}_\theta^i}{R^2 \partial \theta^2} + \frac{1-v}{2} \left( -\frac{\partial^2 \bar{\psi}_r^i}{R \partial R \partial \theta} + \frac{\partial \bar{\psi}_r^i}{R^2 \partial \theta} - \frac{\bar{\psi}_\theta^i}{R^2} + \frac{\partial \bar{\psi}_\theta^i}{R \partial R} + \frac{\partial^2 \bar{\psi}_\theta^i}{\partial R^2} \right) \right]$$

$$-S_2^i (\bar{\psi}_\theta^i + \frac{\partial \bar{w}^i}{R \partial \theta}) = -S_1^i \lambda_i^2 (\bar{I}_2^i \bar{v}_0^i + \bar{I}_3^i \bar{\psi}_\theta^i)$$

(9d)

$$\delta_i^2 S_2^i \left[ \left( \frac{\partial \bar{\psi}_r^i}{\partial R} + \frac{\bar{\psi}_r^i}{R} + \frac{\partial \bar{\psi}_\theta^i}{R \partial \theta} \right) + \bar{\Delta} \bar{w}^i \right] = -S_1^i \lambda_i^2 \bar{I}_1^i \bar{w}^i$$

(9e)

where

$$(\bar{I}_1^i, \bar{I}_2^i, \bar{I}_3^i) = \int_{-1/2}^{1/2} \frac{\rho(z)}{\rho_c} (1, Z_i, Z_i^2) dZ_i$$

(10)

$$(\bar{K}_1^i, \bar{K}_2^i, \bar{K}_3^i) = \frac{1}{E_c h_i^k} \int_{-h_i/2}^{h_i/2} \frac{E(z)}{1-v^2} (1, z, z^2) dz \quad , \quad k = 1, 2, 3$$

(11)

$$\bar{\Delta} = \frac{\partial^2}{\partial R^2} + \frac{\partial}{R \partial R} + \frac{\partial^2}{R^2 \partial \theta^2}$$

(12)

$$S_1^i = \frac{\delta_i^2}{12(1-v^2)} \quad , \quad S_2^i = \frac{\kappa^2 (1-v) \bar{K}_1^i}{2\delta_i^2} \quad , \quad \lambda_i = \omega r_n^2 \sqrt{\frac{\rho_c h_i}{D_i}} \quad , \quad D_i = \frac{E_c h_i^3}{12(1-v^2)}$$

(13a-d)

In order to solve displacement field $\bar{w}^i$, two auxiliary functions are defined as follows

$$\Phi_1^i = \frac{\partial \bar{u}_0^i}{\partial R} + \frac{\bar{u}_0^i}{R} + \frac{\partial \bar{v}_0^i}{R \partial \theta} \quad , \quad \Phi_2^i = \frac{\partial \bar{\psi}_r^i}{\partial R} + \frac{\bar{\psi}_r^i}{R} + \frac{\partial \bar{\psi}_\theta^i}{R \partial \theta}$$

(14a, b)

Using these Eqs. (14a, b) and after some mathematical manipulation, finally a sixth order partial differential equation with constant coefficients is acquired in terms of $\bar{w}^i$ as follows

$$A_1^i \bar{\Delta}\bar{\Delta}\bar{\Delta} \bar{w}^i + A_2^i \bar{\Delta}\bar{\Delta} \bar{w}^i + A_3^i \bar{\Delta} \bar{w}^i + A_4^i \bar{w}^i = 0$$

(15)

where the coefficients $A_k^i$ (k=1,2,3,4) are determined by

$$A_1^i = S_2^i \delta_i^2 \left( \bar{K}_1^i \bar{K}_3^i - \bar{K}_2^{i\,2} \right)$$

(16a)

$$A_2^i = S_1^i \lambda_i^2 \left[ \bar{I}_1^i \left( \bar{K}_1^i \bar{K}_3^i - \bar{K}_2^{i\,2} \right) + S_2^i \delta_i^2 \left( \bar{K}_1^i \bar{I}_3^i + \bar{K}_3^i \bar{I}_1^i - 2\bar{K}_2^i \bar{I}_2^i \right) \right] \tag{16b}$$

$$A_3^i = S_1^i \lambda_i^2 \left[ S_1^i \lambda_i^2 \left( \bar{I}_1^i \left( \bar{K}_1^i \bar{I}_3^i + \bar{K}_3^i \bar{I}_1^i - 2\bar{K}_2^i \bar{I}_2^i \right) + S_2^i \delta_i^2 \left( \bar{I}_1^i \bar{I}_3^i - \bar{I}_2^{i\,2} \right) \right) - S_2^i \bar{K}_1^i \bar{I}_1^i \right] \tag{16c}$$

$$A_4^i = S_1^{i\,2} \lambda_i^4 \bar{I}_1^i [S_1^i \lambda_i^2 \left( \bar{I}_1^i \bar{I}_3^i - \bar{I}_2^{i\,2} \right) - S_2^i \bar{I}_1^i ] \tag{16d}$$

The function $\bar{w}^i(R,\theta)$ can be written as

$$\bar{w}^i(R,\theta) = \hat{w}^i(R)\cos(p\theta) \tag{17}$$

in which the non-negative integer $p$ represents the circumferential wave number of the corresponding mode shape. By substituting Eq. (17) into Eq. (15), the reduced following form of the equation is obtained as

$$\left(\hat{\Delta} - x_1^i\right)\left(\hat{\Delta} - x_2^i\right)\left(\hat{\Delta} - x_3^i\right)\hat{w}^i(R) = 0 \tag{18}$$

where $\hat{\Delta} = \dfrac{\partial^2}{\partial R^2} + \dfrac{\partial}{R\partial R} - \dfrac{p^2}{R^2}$ and $x_1^i$, $x_2^i$, and $x_3^i$ are the roots of the following equation

$$A_1^i x^3 + A_2^i x^2 + A_3^i x + A_4^i = 0 \tag{19}$$

After solving the obtained third order Eq. (19), finally the general solution of Eq. (15) can be expressed as the summation of three Bessel functions as follows

$$\bar{w}^i(R,\theta) = \sum_{k=1}^{3} [c_k^i w_{k1}^i(p, \chi_k^i R) + c_{k+3}^i w_{k2}^i(p, \chi_k^i R)]\cos(p\theta) \tag{20}$$

where $\chi_k^i = \sqrt{|x_k^i|}$, and

$$w_{k1}^i = \begin{cases} J_p(\chi_k^i R), & x_k^i < 0 \\ I_p(\chi_k^i R), & x_k^i > 0 \end{cases}, \quad w_{k2}^i = \begin{cases} Y_p(\chi_k^i R), & x_k^i < 0 \\ K_p(\chi_k^i R), & x_k^i > 0 \end{cases} \quad k=1,2,3 \tag{21a,b}$$

$c_k^i$ are unknown coefficients and $J_p$ and $Y_p$ are the Bessel functions of the first and second kind, respectively, whereas $I_p$ and $K_p$ are the modified Bessel functions of the first and second kind, respectively.

In order to solve displacement fields $\bar{u}_0^i$, $\bar{v}_0^i$, $\bar{\psi}_r^i$ and $\bar{\psi}_\theta^i$, four auxiliary functions $f_1^i$, $f_2^i$, $f_3^i$ and $f_4^i$ [7] are introduced as follows

$$\begin{aligned}
f_1^i &= \bar{K}_1^i \bar{u}_0^i + \bar{k}_2^i \bar{\psi}_r^i & f_2^i &= \bar{K}_1^i \bar{v}_0^i + \bar{K}_2^i \bar{\psi}_\theta^i \\
f_3^i &= \bar{K}_2^i \bar{u}_0^i + \bar{k}_3^i \bar{\psi}_r^i & f_4^i &= \bar{K}_2^i \bar{v}_0^i + \bar{K}_3^i \bar{\psi}_\theta^i
\end{aligned} \tag{22a-d}$$

In order to determine $f_1^i$, $f_2^i$, $f_3^i$ and $f_4^i$, the following forms of solution are considered

$$f_1^i = a_1^i \frac{\partial \bar{w}_1^i}{\partial R} + a_2^i \frac{\partial \bar{w}_2^i}{\partial R} + a_3^i \frac{\partial \bar{w}_3^i}{\partial R} + a_4^i \frac{\partial \bar{w}_4^i}{R\partial\theta} + a_5^i \frac{\partial \bar{w}_5^i}{R\partial\theta} \tag{23a}$$

$$f_2^i = b_1^i \frac{\partial \bar{w}_1^i}{R\partial\theta} + b_2^i \frac{\partial \bar{w}_2^i}{R\partial\theta} + b_3^i \frac{\partial \bar{w}_3^i}{R\partial\theta} + b_4^i \frac{\partial \bar{w}_4^i}{\partial R} + b_5^i \frac{\partial \bar{w}_5^i}{\partial R} \tag{23b}$$

$$f_3^i = a_6^i \frac{\partial \bar{w}_1^i}{\partial R} + a_7^i \frac{\partial \bar{w}_2^i}{\partial R} + a_8^i \frac{\partial \bar{w}_3^i}{\partial R} + a_9^i \frac{\partial \bar{w}_4^i}{R\partial \theta} + a_{10}^i \frac{\partial \bar{w}_5^i}{R\partial \theta} \tag{23c}$$

$$f_4^i = b_6^i \frac{\partial \bar{w}_1^i}{R\partial \theta} + b_7^i \frac{\partial \bar{w}_2^i}{R\partial \theta} + b_8^i \frac{\partial \bar{w}_3^i}{R\partial \theta} + b_9^i \frac{\partial \bar{w}_4^i}{\partial R} + b_{10}^i \frac{\partial \bar{w}_5^i}{\partial R} \tag{23d}$$

in which $a_k^i$ and $b_k^i$ are unknown coefficients. Also, $\bar{w}_4^i$ and $\bar{w}_5^i$ are unknown functions. By substituting Eqs. (22a-d) into Eqs. (9a-e) and using the solutions (23a-d), the coefficients $a_k^i$ and $b_k^i$ as well as the functions $\bar{w}_4^i$ and $\bar{w}_5^i$ can be determined as follows

$$a_k^i = b_k^i = \begin{cases} \dfrac{G_2^i G_5^i}{x_k^{i\,2} - (G_1^i + G_4^i)x_k^i + G_1^i G_4^i - G_2^i G_3^i}, & k = 1,2,3 \\ \dfrac{(x_{k-5}^i - G_1^i)a_{k-5}^i}{G_2^i}, & k = 6,7,8 \end{cases} \tag{24}$$

where

$$G_1^i = \frac{S_1^i \lambda_i^2 (\bar{K}_3^i \bar{I}_1^i - \bar{K}_2^i \bar{I}_2^i)}{\bar{K}_2^{i\,2} - \bar{K}_1^i \bar{K}_3^i} \qquad G_2^i = \frac{S_1^i \lambda_i^2 (\bar{K}_1^i \bar{I}_2^i - \bar{K}_2^i \bar{I}_1^i)}{\bar{K}_2^{i\,2} - \bar{K}_1^i \bar{K}_3^i}$$

$$G_3^i = \frac{S_1^i \lambda_i^2 (\bar{K}_3^i \bar{I}_2^i - \bar{K}_2^i \bar{I}_3^i) + S_2^i \bar{K}_2^i}{\bar{K}_2^{i\,2} - \bar{K}_1^i \bar{K}_3^i} \qquad G_4^i = \frac{S_1^i \lambda_i^2 (\bar{K}_1^i \bar{I}_3^i - \bar{K}_2^i \bar{I}_2^i) - S_2^i \bar{K}_1^i}{\bar{K}_2^{i\,2} - \bar{K}_1^i \bar{K}_3^i}$$

$$G_5^i = S_2^i \qquad G_6^i = \frac{S_2^i \delta_i^2 \bar{K}_2^i}{\bar{K}_2^i - \bar{K}_1^i \bar{K}_3^i}$$

$$G_7^i = \frac{-S_2^i \delta_i^2 \bar{K}_1^i}{\bar{K}_2^i - \bar{K}_1^i \bar{K}_3^i} \qquad G_8^i = -S_1^i \lambda_i^2 \bar{I}_1^i \tag{25a-h}$$

$$\bar{w}_k^i = \hat{w}_k^i \sin(p\theta), \qquad k = 4,5 \tag{26}$$

$$\hat{w}_4^i = c_7^i w_{41}^i(p, \chi_4^i R) + c_8^i w_{42}^i(p, \chi_4^i R) \tag{27a}$$

$$\hat{w}_5^i = c_9^i w_{51}^i(p, \chi_5^i R) + c_{10}^i w_{52}^i(p, \chi_5^i R) \tag{27b}$$

in which $\chi_k^i = \sqrt{|x_k^i|}$, and

$$x_4^i = \frac{\xi_1^i + \sqrt{\xi_1^{i\,2} - 4\xi_2^i}}{2}, \qquad x_5^i = \frac{\xi_1^i - \sqrt{\xi_1^{i\,2} - 4\xi_2^i}}{2} \tag{28a,b}$$

$$\xi_1^i = \frac{2(G_1^i + G_4^i)}{1-\nu}, \qquad \xi_2^i = \frac{4(G_1^i G_4^i - G_2^i G_3^i)}{(1-\nu^2)} \tag{29a,b}$$

$$w_{k1}^i = \begin{cases} J_p(\chi_k^i R), & x_k^i < 0 \\ I_p(\chi_k^i R), & x_k^i > 0 \end{cases}, \qquad w_{k2}^i = \begin{cases} Y_p(\chi_k^i R), & x_k^i < 0 \\ K_p(\chi_k^i R), & x_k^i > 0 \end{cases} \qquad k = 4,5 \tag{30a,b}$$

$$a_k^i = -b_k^i = \begin{cases} \dfrac{G_2^i}{(\dfrac{1-v}{2})x_k^i - G_1^i}, & k = 4,5 \\ 1, & k = 9,10 \end{cases} \quad (31)$$

Finally, the exact solutions for $\bar{u}_0^i$, $\bar{v}_0^i$, $\bar{\psi}_r^i$ and $\bar{\psi}_\theta^i$ according to Mindlin's theory, are obtained by using Eqs. (22a-d). In order to satisfy the continuity conditions, the dimensional forms of the displacement components and stress resultants must be used instead of the dimensionless ones. This is because the annular segments of the stepped plate have different values of the thicknesses $h_i$. Based on the FSDT, ten continuity conditions should be satisfied for each step location, which can be written as follows

$$\begin{aligned} w^i = w^{i+1}, \quad u_0^i = u_0^{i+1}, \quad v_0^i = v_0^{i+1}, \quad \psi_r^i = \psi_r^{i+1}, \quad \psi_\theta^i = \psi_\theta^{i+1} \\ Q_r^i = Q_r^{i+1}, \quad N_r^i = N_r^{i+1}, \quad N_\theta^i = N_\theta^{i+1}, \quad M_r^i = M_r^{i+1}, \quad M_\theta^i = M_\theta^{i+1} \end{aligned} \quad (32)$$

Both inner and outer edges of the stepped circular/annular plate can take any combinations of classical boundary conditions, including free, soft simply supported, hard simply supported and clamped. It should be noted that, for stepped circular FG plates, second type of Bessel function becomes singular at $r=0$. Hence, the unknown coefficients of them ($c_4^i, c_5^i, c_6^i, c_8^i, c_{10}^i$) should be equal to zero. The satisfaction of the continuity and boundary conditions lead to a coefficient matrix. For a nontrivial solution, the determinant of the coefficient matrix must be set to zero for each $p$ and solving this obtained eigenvalue equations yields the frequency parameters $\beta$.

## 3. Results and Discussion

This section contains two parts; firstly, the authors try to validate the present solution with finite element analysis. After verification of results, the effects of the geometrical properties of stepped thickness circular/annular plates on the frequency parameters will be discussed. It should be noted that in this paper Poisson's ratio is assumed to be 0.3 and the shear correction factor has been taken to be 5/6. The numbers in parentheses *(m, n)* show that the vibrating mode has *m* nodal diameters and vibrates in the *n*th mode for the given *m* value. The type of FG plate which is used in this study has the Young's modulus and mass density $E_m = 70 GPa$, $\rho_m = 2700 kg/m^3$ and $E_c = 380 GPa$, $\rho_c = 3800 kg/m^3$ for metallic and ceramic constituents, respectively and power-law index *g* sets to be 1 in all calculations. Also, in this paper, the dimensionless frequency parameter $\beta = \omega r_n^2 \sqrt{\rho_c h_n / D_n}$ is defined in order to evaluate the effects of geometrical properties of stepped thickness circular plates on the natural frequencies.

A comparative study for evaluation of first ten natural frequencies (Hz) of stepped FG circular /annular plate between the present exact solution and the finite element analysis is carried out in

**Table 1.** Numerical results have been calculated for free, soft-simply supported and clamped FG circular Mindlin plates with one step variation. The step locations and thicknesses are selected as $h_1 = 0.2(m)$, $r_1 = 1(m)$, $h_2 = 0.1(m)$, and $r_2 = 2(m)$. To demonstrate further the high accuracy of the present exact solution, in **Table 1**, ANSYS software package of version 12 and ABAQUS software package of version 6.10 are used to model stepped FG circular plate. The plate is analysed with a shell element of type Shell 281 in ANSYS and a solid element in ABAQUS, created on the basis of the FSDT and the 3D elasticity, respectively. A mesh sensitivity analysis was carried out to ensure independency of finite element (FE) results from the number of elements. Results of **Table 1** reveal that very good agreement is achieved for the circular FG plate.

**Table 1.** Comparison study of first ten natural frequencies (Hz) for FG circular plate with one step variation

| B. C's | Mode number (m,n) | (2,1) | (0,1) | (3,1) | (4,1) | (1,1) | (5,1) | (2,2) | (6,1) | (0,2) | (7,1) |
|---|---|---|---|---|---|---|---|---|---|---|---|
| Free | Present | 83.543 | 128.289 | 146.398 | 227.583 | 230.416 | 331.664 | 395.639 | 457.797 | 477.222 | 604.337 |
|  | FEM (FSDT) | 83.578 | 128.340 | 146.460 | 227.670 | 230.510 | 331.800 | 395.790 | 458.000 | 477.400 | 604.660 |
|  | FEM (3D) | 82.773 | 125.658 | 145.156 | 225.139 | 226.525 | 331.018 | 392.833 | 457.637 | 477.737 | 604.696 |
|  | Mode number (m,n) | (0,1) | (1,1) | (2,1) | (0,2) | (3,1) | (1,2) | (4,1) | (5,1) | (2,2) | (0,3) |
| Soft-Simply Supported | Present | 58.084 | 144.063 | 297.991 | 382.925 | 468.489 | 606.805 | 628.476 | 790.895 | 796.988 | 867.957 |
|  | FEM (FSDT) | 58.108 | 144.120 | 298.110 | 383.070 | 468.650 | 607.000 | 628.690 | 791.170 | 797.230 | 868.220 |
|  | FEM (3D) | 57.317 | 142.386 | 297.923 | 383.924 | 469.695 | 603.836 | 630.583 | 790.894 | 794.823 | 860.899 |
|  | Mode number (m,n) | (0,1) | (1,1) | (2,1) | (0,2) | (3,1) | (1,2) | (4,1) | (5,1) | (2,2) | (0,3) |
| Clamped | Present | 110.629 | 223.279 | 392.735 | 493.537 | 594.290 | 767.116 | 782.061 | 958.938 | 978.960 | 1044.083 |
|  | FEM (FSDT) | 110.670 | 223.360 | 392.870 | 493.710 | 594.480 | 767.330 | 782.290 | 959.220 | 979.210 | 1044.300 |
|  | FEM (3D) | 109.880 | 220.465 | 392.132 | 495.484 | 596.716 | 767.374 | 785.848 | 963.867 | 972.748 | 1033.770 |

**Figure 2** shows the variations of the first three frequency parameters $\beta$ versus the step location $\Re = r_1/r_2$ for free circular FG Mindlin plate ($g=1$) with one step variation. The step thickness ratio $h_1/h_2$ and plate thickness ratio $h_2/r_2$ are set to be 3/2 and 0.1, respectively. **Figure 2** shows that as the increase of the step location $\Re$, the first three frequency parameters of free stepped circular FG plate increase slowly to their maximum values and then decrease to the values corresponding to a free circular FG plate without step variation. The maximum values of the first, second and third modes are around *0.8*, *0.85* and *0.9*, respectively.

**Figure 3** shows the variations of the first three frequency parameters $\beta$ of simply supported circular FG plate with one step variation versus the step thickness ratio $\tau = h_1/h_2$. The step location $r_1/r_2$ and the plate thickness ratio $h_2/r_2$ are set to 0.5 and 0.05, respectively. As it is shown in **Figure 3**, For third mode of the simply supported FG plate, as the step thickness ratio $\tau$ enhances, the frequency parameter $\beta$ increases, keeping all other parameters fixed. But the frequency parameters of the first and second modes of the simply supported plate increase slowly to their maximum val-ues and then decrease, indicating that the plate stiffness corresponding to the first and second modes are maximum at the points around $\tau = 2.2$ and $\tau = 1.8$, respectively.

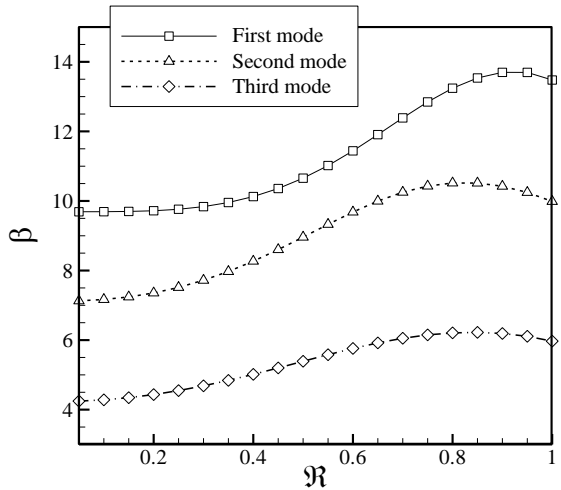
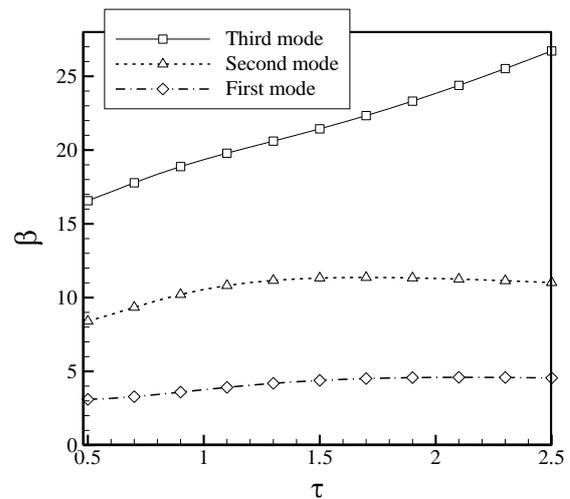

**Figure 2**. Variation of the first three frequency parameters $\beta$ versus the step location $\Re = r_1/r_2$ for free circular FG Mindlin plate with one step variation.

**Figure 3**. Variation of the first three frequency parameter $\beta$ versus the step location $\tau = h_1/h_2$ for simply supported circular FG plate with one step variation.

## 4. Conclusion

The main objective of this paper was to develop an analytical procedure for solving the free vibration of stepped circular FG plates. By introducing some auxiliary and potential functions, the domain decomposition method is employed to find the free dynamic response of stepped FG plate. The proposed model is validated by comparing its dynamic results with finite element analysis and finally, the influence of different parameters of the stepped FG plate such as step locations and step thickness ratios on the natural frequencies is investigated. Results show that the step parameters play a significant role in the determination of vibration behavior of the FG plate, especially for higher vibrating modes.